\documentclass[12pt,leqno,twoside]{article}

\usepackage{bbm}
\usepackage{amssymb}
\usepackage{amsthm}
\usepackage{amsfonts}
\usepackage{graphicx}
\usepackage{verbatim}
\usepackage{enumerate}
\usepackage[intlimits]{amsmath}
\usepackage{fullpage}

\usepackage[english]{babel}

\newtheorem{theorem}{Theorem}
\numberwithin{theorem}{section}
\newtheorem{corollary}{Corollary}
\numberwithin{corollary}{section}
\newtheorem{lemma}{Lemma}
\numberwithin{lemma}{section}

\newtheorem{definition}{Definition}
\numberwithin{definition}{section}
\theoremstyle{definition}
\newtheorem{example}{Example}
\numberwithin{example}{section}

\newcommand{\RR}{\mathbb{R}}
\newcommand{\na}{\mathbf{a}}
\newcommand{\nb}{\mathbf{b}}

\newcommand{\Sym}{\mathrm{Sym}}
\newcommand{\Pd}{\mathrm{Pd}}
\newcommand{\du}{\mathrm{d}}
\newcommand{\Ec}{\mathcal{E}}

\allowdisplaybreaks

\begin{document}

\title{Lower threshold ground state energy 
and testability of minimal balanced cut density}

\author{Andr\'as Kr\'amli \thanks{Bolyai Institute, University of Szeged. The Project is supported by the Hungarian Scientific Research Fund OTKA 105645. E-mail: \textrm{kramli@math.u-szeged.hu}}
\and 
Roland Mark\'o\thanks{Hausdorff Center for Mathematics, University of Bonn. Supported in part by a Hausdorff scholarship. E-mail: \textrm{roland.marko@hcm.uni-bonn.de}}}


\date{}

\maketitle

\begin{abstract} 
Lov\'asz and his coauthors in \cite{lovlac2}  
defined the  notion of microcanonical ground 
state energy $\hat{\mathcal{E}}_\na (G,J)$ -- borrowed from the statistical physics -- for 
weighted graphs $G$, where $\na\in \Pd_q$ is a probability 
distribution on $\{1,...,q\}$  and $J$ is a symmetric $q \times q$ matrix with real entries.
We define a new version of the ground state energy $\hat{\mathcal{E}}^c (G,J)=\inf_{\na\in A_c}\hat{\mathcal{E}}_\na (G,J)$,
called lower threshold ground state energy,  where $A_c = \{\na\in \Pd_q :\, a_i\ge c,\,i=1,\dots, q \}$. Both types of 
energies can be extended for graphons $W$, the limit objects of convergent sequences of simple graphs. The main result 
of the paper 
is Theorem \ref{implik} stating that if $0\leq c_1<c_2 \leq 1$, then the convergence of the sequences $(\hat{\mathcal{E}}^{c_2/q} (G_n,J))$ for each $J \in \Sym_q$
implies convergence of the sequences  $(\hat{\mathcal{E}}^{c_1/q} (G_n,J))$ for each $J \in \Sym_q$.
 As a byproduct one can derive in a natural way the testability of minimum balanced multiway cut densities --  one
 of the fundamental problems of cluster analysis -- proved in 
\cite{BKK}.
\end{abstract}

\section{Preliminaries and notation}

The main goal of this paper is to introduce and to reveal the properties of an intermediate object  between the 
microcanonical ground state energy (MGSE) and ground state energy (GSE) of weighted graphs defined in \cite{lovlac2}.  
For this purpose we need to define the fundamental notions used in
 \cite{lovlac},\cite{lovlac2}, and \cite{lovszeg} and to cite the 
main results therein necessary for us, these will be given below.
 The main contribution of the paper is that we give a
 convergence hierarchy with respect to the aforementioned intermediate
 objects that are Hamiltonians subject to certain conditions. 
This can be regarded as a refined version of Theorem 2.9. (ii) from \cite{lovlac2}
 combined with the equivalence assertion of
 Theorem 2.8. (v) from the same paper. In short,
 these state that MGSE convergence implies GSE convergence.
 Counterexamples are provided indicating that the implication is strict.
 We also reprove with the aid of the established hierarchy
one
of the main results of $\cite{BKK}$. Our motivation comes from 
cluster analysis, where the minimal cut problem is a central
 subject of research. The graph limit theory of Lov\'asz et al. sheds 
new light on this, especially their statistical physics correspondence suits for application in the cluster analysis setting. 

We consider both unweighted simple graphs $G$ (graphs without loops  and multiple edges) and  weighted graphs. We denote 
 the  node and edge sets of $G$ by $V(G)$ and  $E(G)$, respectively. Usually we denote by $\alpha_i= \alpha_i(G)>0$ the weight
  associated with the node $i$  and  $\beta_{ij}=\beta_{ij} (G)\in \RR $ the weight associated with the edge $ij$. 
 We set $\alpha_G=\sum_i\alpha_i(G)$.

Definitions 1.1. -- 1.9. (except  1.7.) are borrowed from  \cite{lovlac} and \cite{lovlac2}. In order to facilitate the reading, at every definition we indicate its exact place of  occurrence. 
\begin{definition}\label{grafon}(\cite{lovlac} Definition 3.1.) 
Let $\mathcal{W}$ denote the space  of bounded symmetric measurable functions $W\colon [0,1]^2 \to \RR$, that is  $W(x,y)=W(y,x)$. 
 Assume that the functions $W\in \mathcal{W}$ take their values in an interval $I$, usually $I=[0,1]$.
We can think of the interval $[0,1]$ as the set of nodes of graph with a node set that has  cardinality continuum, and of the values $W(x,y)$ as the weight of the edge $xy$.
We call the functions in $\mathcal{W}_I$ {\it{graphons}.}
\end{definition}

\begin{definition}\label{cutdist} (\cite{lovlac2} Definition 2.4.)
Let $G$ and $G'$ be two weighted graphs with node set $V$ and $V'$, respectively.
For $i\in V$  and $u \in V'$ set $\mu_i=\alpha_i(G)/\alpha_G$ and  ${\mu'}_u=\alpha_u(G')/\alpha_{G'}$. Then we define the
set of fractional overlays $\chi(G,G')$
as the set
 of probability
 distributions $X$ on $V\times V'$ (or couplings of $\mu$ and $\mu'$) such that
$$
\sum_{u\in V'} X_{iu}=\mu_i \textrm{ for all } i\in V, \textrm{ and }
 \sum_{i\in V}X_{iu}={\mu'}_u \textrm{ for all } u\in V', 
$$

and set
\begin{equation}\label{kockadist}
\delta_\square(G,G')=\min_{X\in \chi(G,G') }\max_{S, T\subset V\times V'}
\left\vert \sum_{\substack{(i,u)\in S \\(j,v)\in T}}  X_{iu} X_{jv} \left(\beta_{ij}(G)       
- \beta_{uv}(G')\right) \right\vert.
\end{equation}

\end{definition}

\begin{definition}\label{cutnormgrafon} (\cite{lovlac} formula (3.3))
The distance defined in  Definition \ref{cutdist} can be extended to graphons  in terms of the {\it{cut norm}}
\begin{align}\label{cutnorm}
\Vert W \Vert_\square &= \sup_{S, T\subset [0,1]}\left\vert \int_{S\times T}W(x,y)dxdy \right \vert\\
&=\sup_{f,g\colon [0,1]\to [0,1]}\left\vert \int W(x,y)f(x)g(y)dxdy\right \vert, \nonumber
\end{align}

where the suprema go over measurable subsets and functions, respectively.
\end{definition}

\begin{definition}\label{cutdistgrafon}(\cite{lovlac2} formula (3.3)) 
The cut distance of two graphons $U$ and $W$
is defined as 
\begin{equation}\label{cutdistgraf}
\delta_\square(U,W)=\inf_\phi \Vert U-W^\phi \Vert_\square,
\end{equation}where the infimum goes over all measure preserving permutations of $[0,1]$, $W^\phi(x,y)=W(\phi(x),\phi(y))$.
\end{definition}

Now we define
 three versions of the {\it{ground state energies}} (GSE) borrowed from the statistical physics 
that are the objects investigated in this 
paper (for the mathematical treatment of statistical physics, see, e.g. Sinai's book \cite{S}). 
They are defined  in terms of a finite set of states $[q]=\{1,\dots ,q\}$,
 and a symmetric $q\times q$ matrix $J$ with entries 
in $\RR$, the set of these matrices is denoted by $\textrm{Sym}_q$. 
A {\it{spin configuration}} on a simple or weighted graph $G$
is given by a map 
$\phi \colon V(G) \to [q]$.
\begin{definition}\label{gse} (\cite{lovlac2} formulae (2.8) and (2.11))

\noindent {\it{The energy density of a spin configuration of $G$ with respect to $J$ is given by}}
\begin{equation}\label{energiasur}
\Ec_\phi(G,J)=-\frac{2}{\left|V(G)\right|^2} \sum_{uv\in E(G)}J_{\phi(u)\phi(v)}.
\end{equation}

\noindent {\it{The ground state energy (GSE) of $G$ with respect to $J$ is }} 
\[
\hat{\Ec}(G,J)= \min_{\phi\colon V(G)\to [q]}{\Ec_\phi}(G,J).
\]
For a graphon $W$ and a $\rho=(\rho_1, \dots, \rho_q)$ with $\rho_i\colon [0,1] \to [0,1]$ being measurable and satisfying $\sum_i \rho(x)=1$ for each $x \in [0,1]$ called $q$-fractional partition the energy is defined as
$$
\Ec_\rho(W,J)=-\sum_{i,j=1}^{q} J_{ij} \int_{[0,1]^2} \rho_i(x) \rho_j(y) W(x,y) \du x \du y,
$$ and the GSE is 
$$
\Ec(W,J)= \min_{\rho}{\Ec_\rho}(W,J).
$$   
\end{definition}

Let $\Pd_q$ be the set of all probability distributions on $[q]$.

Imposing some restrictions on the set where the minimum is taken in the above definition we can define another version of energies that are important in graph limit theory.
\begin{definition}\label{mgse}(\cite{lovlac2} formula (2.14))
{\it{Microcanonical ground state energy (MGSE) of $G$ with respect to $J \in \Sym_q$ and a probability distribution $\na=(a_1,\dots,a_q) \in \Pd_q$ is defined using the set}} \\

\begin{equation}\label{afelosztas}
\Omega_\na(G)=
\left\{
\phi\colon V(G) \to [q]\colon  \left|  |\phi^{-1}(\{i\})| -a_i |V(G)|  \right| \leq 1 
\textrm{ for  all  } i\in [q] \right\},
\end{equation}
and is the quantity 
\begin{equation}\label{mikrokan}
\hat{\Ec}_\na(G,J) =\min_{\phi \in \Omega_\na(G)}\Ec_\phi(G,J).
\end{equation}
Let $\omega_\na=\left\{ \rho: \int_0^1 \rho_i(x)\du x = a_i \textrm{ for  all  } i\in [q] \right\}$ be a subset of $q$-fractional partitions. Then the MGSE of a graphon $W$ is defined as
$$
\Ec_\na(W,J) =\min_{\rho \in \omega_\na}\Ec_\rho(W,J).
$$
\end{definition}

Next we introduce the central object of our current investigation.
\begin{definition}\label{ltgs}
Let $G$ be a weighted graph, $q \geq 1$, $J\in \Sym_q$, and $0 \leq c \leq 1/q$. 
 We define the set $A_c = \{\na \in \Pd_q :\quad a_i \geq c, i=1, \ldots , q\}$, and with its help the lower threshold 
ground state energy (LTGSE):
\begin{equation}\label{defltgse1}
\hat{\mathcal{E}}^c(G,J) = \inf_{\na \in A_c} \hat{\mathcal{E}}_\na(G,J).
\end{equation}
In a similar manner we introduce the LTGSEs for a graphon $W$ for $q \geq 1$, $J\in \Sym_q$, and $0 \leq c \leq 1/q$ 
lower threshold:
\begin{equation}\label{defltgse}
\mathcal{E}^c(W,J)= \inf_{\na \in A_c} \mathcal{E}_\na(W,J).
\end{equation}
\end{definition}   


We remind the reader of the definition of testability of simple graph parameters. Before doing it, we should define the 
randomization procedure for graphs, used here. 

\begin{definition} \label{defrnd}(\cite{lovlac} Introduction of Section  2.5.3.)
For a graph $G$ and a positive integer $k$ let $\mathbb{G}(k,G) $ denote the random induced subgraph $G[S]$
where $S$ is chosen uniformly from all subsets of $V(G)$ of cardinality $k$.
\end{definition}

\begin{definition} \label{deftest}(\cite{lovlac} Definition 2.11)
A real function $f$ 
defined on the set of simple graphs
 is  a testable simple graph parameter,
 if for every $\varepsilon>0$ there exists a $k=k(\varepsilon) \in \mathbb{N}$ such that for every simple graph $G$ 
on at least $k$ vertices
\[P(\left|f(G)-f(\mathbb{G}(k,G))\right|>\varepsilon)<\varepsilon.\]
\end{definition}


This paper is organized as follows. In the second section we prove yet another equivalent condition to left-convergence of a graph sequence relying on a subclass of MGSE, the reasoning will be instrumental for the proof of our main result in the subsequent section. 
In the third section we study the convergence of LTGSEs, see (\ref{defltgse}) for their definition. 
We will consider $c \colon \mathbb{N} \to [0,1]$ threshold functions with the property that  $c(q)q$ is constant as a function of $q$. For this case we will prove that if $0 \leq c_1(q) < c_2(q) \leq 1/q$ (for all $q$), then the convergence of $(\mathcal{E}^{c_2(q)}(W_{n},J))_{n \geq 1}$ for all $q \geq 1$, $J \in \Sym_q$ implies the convergence of $(\mathcal{E}^{c_1(q)}(W_n,J))_{n \geq 1}$ for all $q \geq 1$, $J \in \Sym_q$. 

In the fourth section we provide some examples of graphs and graphons which support the fact, that the implication of convergence in the third section is strict in the sense that convergence of LTGSE sequences with smaller threshold do not imply convergence of  LTGSE sequences with larger threshold in general. We also present a one-parameter family of block-diagonal graphons whose elements can 
be
distinguished by LTGSEs for any threshold $c >0$.
\section{Microcanonical convergence}

We start by showing that for each discrete probability distribution with rational probabilities there exists a uniform probability distribution, such that the microcanonical ground state energies (MGSE) of it can be expressed as MGSEs corresponding to the uniform distribution.
\begin{lemma} \label{racpart}
Let $q \geq 1$, and $\na \in \Pd_{q}$ be  such that $\na=(\frac{k_{1}}{q'},\frac{k_{2}-k_{1}}{q'}, \dots, \frac{k_{q}-k_{q-1}}{q'})$, where $q'$ is a positive integer, $k_{1} \leq k_{2} \leq \dots \leq k_{q}=q'$ are non-negative integers. Then for all $J \in \Sym_{q}$ there exists a $J' \in \Sym_{q'}$, such that for $\nb=(1/q', \dots , 1/q') \in \Pd_{q'}$ and every graphon $W$ it holds that
\[ \mathcal{E}_{\na}(W,J) = \mathcal{E}_{\nb} (W,J'). \]
\end{lemma}
\proof
Set $k_{0}=0$. If the $i$'th component of $\na$ is $0$, then erase this component from $\na$, and also erase the $i$'th row and column of $J$. This transformation clearly will have no effect on the value of the GSE. 
Let us define the $q' \times q'$ matrix $J'$ by blowing up rows and columns of $J$ in the following way. For each $u, v \in [q']$ let $J'_{uv}= J_{ij}$, where $k_{i-1} < u \leq  k_{i}$ and $k_{j-1} < v  \leq k_{j}$.
The matrix $J'$ defined this way is clearly symmetric.

Now we will show that for every  $q$-fractional partition with distribution $\na$  there exists a $q'$-fractional partition $\rho '$ with distribution $\nb$, and vice versa, such that $\mathcal{E}_{\rho}(W,J) = \mathcal{E}_{\rho '} (W,J')$. On one hand, for $1 \leq u \leq q'$ let $\rho'_{u}= \frac{\rho_{i}}{k_{i}-k_{i-1}}$, where $k_{i-1} < u \leq  k_{i}$. Then

\begin{align*}
\mathcal{E}_{\rho '} (W,J') &= - \sum_{u,v=1}^{q'} J'_{u,v} \int_{[0,1]^2}  \rho'_{u}(x) \rho '_{v}(y) W(x,y) \du x \du y   \nonumber \\
  &= - \sum_{i,j=1}^{q} J_{i,j} \sum_{l=k_{i-1}+1}^{k_{i}} \sum_{h=k_{j-1}+1}^{k_{j}} \int_{[0,1]^2}  \rho '_{l}(x) \rho'_{h}(y) W(x,y) \du x \du y \\ &= 
  \mathcal{E}_{\rho}(W,J).  \nonumber
\end{align*}

On the other hand, for $1 \leq i \leq q$ let $\rho_{i}:= \sum_{l=k_{i-1}+1}^{k{i}} \rho '_{l}$. Then
\begin{align}
 \mathcal{E}_{\rho} (W,J) &= - \sum_{i,j=1}^{q} J_{i,j} \int_{[0,1]^2}  \rho_{i}(x) \rho_{j}(y) W(x,y) \du x \du y  \nonumber \\
 &= - \sum_{i,j=1}^{q} J_{i,j} \sum_{l=k_{i-1}+1}^{k_{i}} \sum_{h=k_{j-1}+1}^{k_{j}} \int_{[0,1]^2}  \rho '_{l}(x) \rho '_{h}(y) W(x,y) \du x \du y  \nonumber \\ 
 &= - \sum_{u,v=1}^{q'} J'_{u,v}  \int_{[0,1]^2}  \rho '_{u}(x) \rho '_{v}(y) W(x,y) \du x \du y = 
 \mathcal{E}_{\rho '}(W,J').  \nonumber
\end{align}

So we conclude that
\[ 
\mathcal{E}_{\na}(W,J) = \inf_{\rho \in \omega_\na} \mathcal{E}_{\rho}(W,J) =  \inf_{\rho ' \in \omega_\nb} \mathcal{E}_{\rho '} (W,J') = \mathcal{E}_{\nb} (W,J').
\]
\qed

First we show that the MGSE with fixed parameters $W$, $J$ are close, whenever their corresponding probability distribution parameters are close to each other.
\begin{lemma}  \label{folyt} 
Let $q \geq 1$, $J \in \Sym_{q}$, and $W$ be an arbitrary graphon. Then for $\na,\nb \in \Pd_{q}$ we have that
\[ 
 \left|\mathcal{E}_{\na}(W,J) - \mathcal{E}_{\nb} (W,J)\right|< 2\|\na-\nb\|_{1} \|W\|_{\infty} \|J\|_{\infty}.
\] 
\end{lemma}
\proof
Let $\rho \in \omega_\na$, 
let us construct according to this a $\rho ' \in \omega_\nb$ the following way. First let us line up those  $i$'s, for which  $b_{i} \geq a_{i}$, for simplicity index them by integers from $1$ to $k$. 
Let $\rho_{1} '$ be such, that $
 \rho_{1} (x) \leq \rho_{1} '(x) \leq 1$ for each $x \in [0,1]$ and $\int_{0}^{1} \rho_{1} '(x) \du x = b_{1}$.
It is clear that such a function exists. We define $ \rho_{2} '$ in similar fashion: let
 $\rho_{2} (x) \leq \rho_{2} '(x) \leq 1 - \rho_{1} '(x)$ for all $x \in [0,1]$ and $\int_{0}^{1} \rho_{2} '(x) \du x = b_{2}$, the existence is again clear. We define subsequently $\rho_{i} '$ for $i$'s obeying $b_{i} \geq a_{i}$ by taking care that $\rho_{i} (x) \leq \rho_{i} '(x) \leq 1 - [\sum_{j=1}^{i-1} \rho_{j} '(x)]$ holds at each step. 
In the other case, when $b_i < a_i$, we reverse the inequality we wish to be satisfied by the functions $\rho_i$ and $\rho'_i$, and define  $\rho_{i} '$ accordingly.
For the constructed $\rho_{i} '$ either $\rho_{i} '(x) \leq \rho_{i}(x)$ for all $x \in [0,1]$, or $\rho_{i} '(x) \geq \rho_{i}(x)$ for all $x \in [0,1]$, and additionally $\sum_i \rho_i(x)=1$. Hence
\begin{align*}
\|\rho - \rho '\|_{1}&= \sum_{i=1}^q \int_0^1 \left|\rho_{i}(x) -\rho '_{i}(x)\right|\du x =
\sum_{i=1}^q \left|\int_0^1 \rho_{i}(x) -\rho '_{i}(x)  \du x\right| \\ &=\sum_{i=1}^q \left|a_{i}-b_{i}\right|= \|\na-\nb\|_{1}.
\end{align*}
Now we give an upper bound on the deviation of MGSEs.
\begin{align*}
&\left|\mathcal{E}_{\rho}(W,J) -\mathcal{E}_{\rho '} (W,J)\right| \\ & \qquad =   \left|\sum_{i,j=1}^q J_{i,j} \int\limits_{[0,1]^2} (\rho_{i}(x) \rho_{j}(y)-\rho_{i} '(x)\rho_{j} '(y))W(x,y) \du x \du y\right|  \nonumber \\
 & \qquad \leq \|W\|_{\infty} \|J\|_{\infty}  \sum_{i,j=1}^q \int\limits_{[0,1]^2} 
\left|\rho_{i}(x) \rho_{j}(y)-\rho_{i}(x) \rho_{j} '(y)\right| \\  & \qquad\quad + \left|\rho_{i} (x)\rho_{j} '(y)- \rho_{i} '(x)\rho_{j} '(y)\right|\du x \du y  \\
  &\qquad \leq \|W\|_{\infty} \|J\|_{\infty}  \sum_{i,j=1}^q a_ {i} \|\rho_{j}-\rho_{j} '\|_1 + b_ {j} \|\rho_{i}-\rho_{i} '\|_1  
\\ & \qquad = 2\|\na-\nb\|_{1} \|W\|_{\infty} \|J\|_{\infty} \nonumber.
\end{align*}
The 
second inequality follows by Fubini's theorem. From the definition of MGSE the statement of the lemma follows.
\qed

With the aid of the two previous lemmas we are able to prove the main assertion of the section.
 In the statement of the following theorem the LTGSE expression $\mathcal{E}^{1/q}(W,J)$ (which is equal to $\mathcal{E}_\nb(W,J)$, 
with $\nb=(1/q, \dots, 1/q)$) appears, the notion will further be generalized in what follows later on.
\begin{theorem}
\label{ekviv}
Let $I$ be a bounded interval, and $(W_{n})_{n \geq 1}$ a sequence of graphons from $\mathcal{W}_I$.
If for all $q \geq 1$ and $J \in \Sym_q$  the sequences $(\mathcal{E}^{1/q}(W_{n},J))_{n \geq 1}$ converge, 
then for all $q \geq 1$, $\na \in \Pd_q$ and $J \in \Sym_q$  the sequences $(\mathcal{E}_{\na}(W_{n},J))_{n \geq 1}$ converge.
\end{theorem}
\proof
Let $q \geq 1$, $\na \in \Pd_q$ and $J \in \Sym_q$ be arbitrary and fixed. We will prove that whenever the conditions of the theorem are satisfied, then $(\mathcal{E}_{\na}(W_{n},J))_{n \geq 1}$ is Cauchy convergent. 
Fix an arbitrary $\varepsilon > 0$. Let $q'$ be such that $ 4\frac{q}{q'}\|I\|_{\infty} \|J\|_{\infty}< \frac{\varepsilon}{3}$, and let $\nb \in \Pd_q$ be such that $b_i=[a_i/q']$ ($i=1, \dots, q-1$), $b_q= 1 -\sum_{i=1}^{q-1} b_i$ (where $[x]$ is the lower integer part $x$). Then 
\[
\|\na-\nb\|_1 = \sum_{i=1}^{q} \left|a_i-b_i\right|\leq 2 \frac{q-1}{q'} <  2 \frac{q}{q'}.
\]
$\nb$ is a $q'$-rational distribution, so by Lemma \ref{racpart} there exists $J' \in \Sym_{q'}$, such that for all $n \geq 1$
\[ 
\mathcal{E}_{\nb}(W_n,J) =\mathcal{E}^{1/q'} (W_n,J').
\]
It follows from the conditions of the theorem that there exists $n_0 \in \mathbb{N}$ such that for all $m,n \geq n_0$ it is true that $\left|\mathcal{E}^{1/q'} (W_n,J')-\mathcal{E}^{1/q'} (W_m,J')\right|< \frac{\varepsilon}{3}$. Applying Lemma \ref{folyt} to all $m,n \geq n_0$ we get that
\begin{align*}
&\left|\mathcal{E}_\na (W_n,J)-\mathcal{E}_\na (W_m,J)\right|\leq 
\left|\mathcal{E}_\na (W_n,J)-\mathcal{E}_\nb (W_n,J)\right|\nonumber \\ 
& \qquad \quad + \left|\mathcal{E}_\nb (W_n,J)-\mathcal{E}_\nb (W_m,J)\right| + \left|\mathcal{E}_\nb (W_m,J)-\mathcal{E}_\na (W_m,J)\right|\nonumber \\
& \qquad \leq 2 \|\na-\nb\|_1 \|I\|_{\infty} \|J\|_{\infty} + \left|\mathcal{E}^{1/q'} (W_n,J')-\mathcal{E}^{1/q'} (W_m,J')\right| \nonumber \\
& \qquad \quad +2 \|\na-\nb\|_1 \|I\|_{\infty} \|J\|_{\infty} \\& \qquad \leq \frac{\varepsilon}{3}+\frac{\varepsilon}{3}+\frac{\varepsilon}{3} =\varepsilon. \nonumber
\end{align*}
\qed

We remark that Theorem \ref{ekviv} also appears in \cite{lovlac2} as Corollary 7.4, but its proof follows a different line of thought in the present paper.

\section{Weaker convergence, lower threshold microcanonical ground state energies}
In various cases of testing, for certain cuts of graphs neither the notion of ground state energies, nor the notion of microcanonical ground state energies are satisfactory. For example when investigating clusteredness of a graph in a certain sense these notions become useless, because the partition for which energies attain the minimal value are trivial partitions. On the other hand, in many applications one only asks for a lower bound on the size of these classes to keep a grade of freedom of the ground state case and  at the same time achieve a certain balance with respect to the sizes of classes. 
This setting can be regarded as an intermediate energy notion that manages to get rid of values corresponding to trivial partitions.
Recall Definition \ref{ltgs} of the lower threshold ground state energies.

The next theorem will deliver an upper bound on the difference of the MGSEs of $G$ and $W_G$ for fixed $\na$ an $J$, $W_G$ is the graphon constructed form the adjacency matrix of $G$ in the natural way. A straightforward consequence of this will be the analogous statement for the LTGSEs.
\begin{theorem}
\cite{lovlac2}
\label{mikrokozel}
Let $G$ be a weighted graph, $q \geq 1$, $\na \in \Pd_q$ and $J \in \Sym_q$. Then 
\[
\left|\hat{\mathcal{E}}_\na(G,J)- \mathcal{E}_\na(W_G,J)\right|\leq 6 q^3 \frac{\alpha_{\max}(G)}{\alpha_G} \beta_{\max} (G) \|J\|_\infty.
\]
\end{theorem}

Since the upper bound in the theorem for a given $q$ is not dependent on  $\na$, it is easily possible to apply it to the LTGSEs.
\begin{corollary}\label{kov}
Let $G$ be a weighted graph, $q \geq 1$, $0 \leq c \leq 1/q$ and $J \in \Sym_q$. Then 
\[
\left|\hat{\mathcal{E}}^c(G,J)- \mathcal{E}^c(W_G,J)\right|\leq 6 q^3 \frac{\alpha_{\max}(G)}{\alpha_G} \beta_{\max} (G) \|J\|_\infty.
\]
\end{corollary}
%

Based on the preceding facts we are able to perform analysis on the LTGSEs the same way as the authors of \cite{lovlac2} did in the case of MGSE.
\begin{corollary}
\label{atmen}
Let $G_n$ be a sequence of weighted graphs with uniformly bounded edge weights. Then if $\frac{\alpha_{\max}(G_n)}{\alpha_{G_n}} \to 0$ ($n \to \infty$), then for all  $q \geq 1$, $0 \leq c \leq 1/q$ and $J \in \Sym_q$ the sequences $(\hat{\mathcal{E}}^c(G_n,J))_{n \geq 1}$ converge if, and only if $(\mathcal{E}^c(W_{G_n},J))_{n \geq 1}$ converge, and then
\[
\lim_{n \to \infty} \hat{\mathcal{E}}^c(G_n,J) =\lim_{n \to \infty} \mathcal{E}^c(W_{G_n},J). 
\]
\end{corollary}

Recall the definition of testability, Definition \ref{deftest}. It was shown in \cite{lovlac}, among presenting other characterizations, that the testability of a graph parameter $f$ is equivalent to the existence of a $\delta_\square$-continuous extension $\hat f$ of $f$ to the space $\mathcal{W}_I$, where extension here means that $f(G_n) - \hat f(W_{G_n}) \to 0$ whenever $\left|V(G_n)\right|\to \infty$ (see \cite{lovlac}, Theorem 6.1, the equivalence of (a) and (d)). Using this we are able to present yet another consequence of Theorem \ref{mikrokozel}, that was verified earlier using a different approach in \cite{BKK} (see also \cite{bolla}, Chapter 4).

\begin{corollary}
For all  $q \geq 1$, $0 \leq c \leq 1/q$ and $J \in \Sym_q$ the simple graph parameter $f(G)=\hat{\mathcal{E}}^c(G,J)$ is testable. Choosing $J$ appropriately, $f(G)$ can be regarded as a type of balanced multiway minimal cut in \cite{BKK}.

\end{corollary}
\proof
Let $q \geq 1$, $0 \leq c \leq 1/q$ and $J \in \Sym_q$ be fixed, and we define $\hat f(W)=\mathcal{E}^c(W,J)$. It follows from Corollary \ref{kov} that $f(G_n) - \hat f(W_{G_n}) \to 0$ whenever $\left|V(G_n)\right|\to \infty$. It remains to show that $\hat f$ is $\delta_\square$-continuous. 
 To elaborate on this issue, let $U,W \in \mathcal W_I$ and $\phi$ be a measure-preserving permutation of $[0,1]$ such that $\delta_\square(U, W)=\|U-W^\phi\|_\square$, and let $\rho=(\rho_1, \dots, \rho_q)$ be an arbitrary fractional partition. Then
\begin{align}
\left|\mathcal{E}_\rho(U,J)-\mathcal{E}_\rho(W^\phi,J)\right| &\leq  \sum_{i,j=1}^q \left|J_{ij}\right| \left|\int_{[0,1]^2}(U-W^\phi)(x,y)\rho_i(x)\rho_j(y) \du x \du y\right| \nonumber \\ 
 & \leq q^2 \|J\|_\infty \|U-W^\phi\|_\square= q^2 \|J\|_\infty \delta_\square(U, W). \label{egy}
\end{align}
This implies our claim, as $\mathcal{E}_\na(W,J)=\mathcal{E}_\na(W^\phi,J)$ for any $\na \in \Pd_q$ and $\phi$ measure preserving permutation, and the fact that the right-hand side of $(\ref{egy})$ does not depend on $\na$, and that by definition $\mathcal{E}^c(W,J)= \inf_{\na \in A_c} \mathcal{E}_\na(W,J)$.

\qed

In order to analyze the convergence relationship of LTGSEs with different thresholds for a given graph sequence it is sensible to consider $c$ as a function of $q$. We restrict our attention to lower threshold functions $c$ with $c(q)q$ being constant, which means that in the case of graphons the total size of the thresholds  stays the same relative to the size of the interval $[0,1]$ (in the case of graphs relative to the cardinality of the vertex set). The main statement of the current section informally asserts that the convergence of LTGSEs with larger lower threshold imply convergence of all LTGSEs with smaller ones.
By the results of the previous section we know that in the case of $c(q)=1/q$ the convergence of these LTGSEs is equivalent convergence of the MGSEs for all probability distributions, and by this, according to \cite{lovlac2}, to left convergence of graphs. Moreover, in the case of $c(q)=0$ it is equivalent to the convergence of the unrestricted GSEs, that property is known to be strictly weaker than left convergence.
For technical purposes we introduce \emph{general} LTGSEs and will refer to the previously presented notion in all that follows as \emph{homogeneous} LTGSEs.
\begin{definition}
Let $q \geq 1$, $\mathbf{x}=(x_1, \dots , x_q)$, $x_1, \dots x_q \geq 0$ and $\sum_{i=1}^q x_i \leq 1$, and let $A_\mathbf{x} = \{\na \in \Pd_q :\quad a_i \geq x_i, i=1, \dots , q\}$. For a graphon $W$ and $J \in \Sym_q$ we call the following expression the lower threshold ground state energy corresponding to $\mathbf{x}$:
\[ \mathcal{E}^\mathbf{x}(W,J)= \inf_{\na \in A_\mathbf{x}} \mathcal{E}_\na (W,J).
\]
The definition of $\hat{\mathcal{E}}^\mathbf{x}(G,J)$ for graphs is analogous.
\end{definition}
Similarly to Lemma \ref{racpart}, the convergence of homogeneous LTGSEs is equivalent to the convergence of certain general LTGSEs. 

\begin{lemma}
\label{ekviv2}
Let $I$ be a bounded interval, $(W_{n})_{n \geq 1}$ a sequence of graphons in $\mathcal{W}_I$.
Let $c$ be a lower threshold function, so that $c(q)q = h$ for all $q \geq 1$ and for some $0 \leq h \leq 1$.
If for all $q \geq 1$ and $J \in \Sym_q$ the sequences $(\mathcal{E}^{c(q)}(W_{n},J))_{n \geq 1}$ converge, then for all $q \geq 1$, all 
\begin{equation} \label{felt}
\mathbf{x}=(x_1, \dots , x_q) \quad x_1, \dots ,x_q \geq 0 \qquad \sum_{i=1}^q x_i = h , 
\end{equation}
 and $J \in \Sym_q$, the sequences $(\mathcal{E}^{\mathbf{x}}(W_{n},J))_{n \geq 1}$ also converge.
\end{lemma}
\proof
Fix an arbitrary graphon $W$ from $\mathcal{W}_I$, $q \geq 1$ and $J \in \Sym_q$, and an arbitrary vector $\mathbf{x}$ that satisfies condition (\ref{felt}). Select for each of these vectors $\mathbf{x}$ a positive vector $\mathbf{x}'$ that obeys the condition (\ref{felt}), and that has components which are integer multiples of $c(q')$ ($q'$ will be chosen later), so that
\[
\|\mathbf{x}-\mathbf{x}'\|_1 \leq 2qc(q') = 2h \frac{q}{q'}.
\]
The sets $A_\mathbf{x}$ and $A_{\mathbf{x}'}$ have Hausdorff distance in the $L^1$-norm at most $\|\mathbf{x}-\mathbf{x}'\|_1$, in particular for every $\na \in A_\mathbf{x}$ there exists a $\nb \in A_{\mathbf{x}'}$, such that $\|\na-\nb\|_1 \leq \|\mathbf{x}-\mathbf{x}'\|_1$, and vice versa. Let $\varepsilon > 0$ be arbitrary, and $\na \in A_\mathbf{x}$ be such that $\mathcal{E}^{\mathbf{x}}(W,J)+\varepsilon >\mathcal{E}_\na(W,J)$ holds. Then by applying Lemma \ref{folyt} we have that
\begin{align*}
\mathcal{E}^{\mathbf{x}'}(W,J)- \mathcal{E}^{\mathbf{x}}(W,J) &<
\mathcal{E}^{\mathbf{x}'}(W,J) -\mathcal{E}_\na(W,J) + \varepsilon \\ &\leq 
\mathcal{E}_{\nb}(W,J) - \mathcal{E}_\na(W,J) + \varepsilon  \\ &\leq
 2\|\na-\nb\|_{1} \|W\|_{\infty} \|J\|_{\infty} + \varepsilon \\ &\leq
2\|\mathbf{x}-\mathbf{x}'\|_{1} \|W\|_{\infty} \|J\|_{\infty} + \varepsilon. \nonumber
\end{align*}
The lower bound of the difference can be handled similarly, and therefore by the arbitrary choice of $\varepsilon$ it holds that 
\[
\left|\mathcal{E}^{\mathbf{x}'}(W,J)- \mathcal{E}^{\mathbf{x}}(W,J)\right|\leq 2\|\mathbf{x}-\mathbf{x}'\|_{1} \|W\|_{\infty} \|J\|_{\infty} \leq 4 h \frac{q}{q'} \|I\|_{\infty} \|J\|_{\infty}.
\]
With completely analogous line of thought to the proof of Lemma \ref{racpart}, one can show that there exists a $J' \in \Sym_{q'}$ such that $\mathcal{E}^{\mathbf{x}'}(W,J) = \mathcal{E}^{c(q')}(W,J')$.
Finally, choose $q'$ small enough in order to satisfy $4 h \frac{q}{q'} \|I\|_{\infty} \|J\|_{\infty} < \frac{\varepsilon}{3}$, and $n_0 > 0$ large enough, so that for all $m,n \geq n_0$ the relation
$$\left|\mathcal{E}^{c(q')}(W_n,J') -\mathcal{E}^{c(q')}(W_m,J')\right|< \frac{\varepsilon}{3}$$ holds.

Then for all $m,n \geq n_0$:

\begin{align*}
&\left|\mathcal{E}^{\mathbf{x}}(W_n,J) -\mathcal{E}^{\mathbf{x}}(W_m,J)\right|<\left|\mathcal{E}^{\mathbf{x}}(W_n,J)- \mathcal{E}^{\mathbf{x}'}(W_n,J)\right| \nonumber \\  & \qquad \quad + \left|\mathcal{E}^{c(q')}(W_n,J') -\mathcal{E}^{c(q')}(W_m,J')\right| 
+ \left|\mathcal{E}^{\mathbf{x}'}(W_m,J)- \mathcal{E}^{\mathbf{x}}(W_m,J)\right| \nonumber \\ & \qquad < \frac{\varepsilon}{3}+\frac{\varepsilon}{3}+\frac{\varepsilon}{3}= \varepsilon. \nonumber
\end{align*}
We did not only prove the statement of the lemma, but we also showed that the convergence is uniform in the sense that $n_0$ does not depend on $\mathbf{x}$ for fixed $q$ and $J$.
\qed

With the aid of the former lemma we can now prove that if all homogeneous LTGSEs with large thresholds converge, then all homogeneous LTGSEs with smaller ones also converge.

\begin{theorem}
\label{implik}
Let $I$ be a bounded interval, $(W_{n})_{n \geq 1}$ a sequence of graphons in $\mathcal{W}_I$.
Let $c_1,c_2$ be two lower threshold functions, so that $c_1(q)q = h_1 < h_2 =c_2(q)q$ for all $q \geq 1$ for some $0 \leq h_1,h_2 \leq 1$.
If for every $q \geq 1$ and $J \in \Sym_q$ the sequences $(\mathcal{E}^{c_2(q)}(W_{n},J))_{n \geq 1}$ converge, then for every $q \geq 1$ and $J \in \Sym_q$ the sequences $(\mathcal{E}^{c_1(q)}(W_{n},J))_{n \geq 1}$ also converge.
\end{theorem}
\proof
From Lemma \ref{ekviv2} it follows that if the conditions of the theorem are satisfied then for every $q \geq 1$, every
\begin{equation} 
\label{feltt}
\mathbf{x}=(x_1, \dots , x_q) \quad x_1, \dots x_q \geq 0 \qquad \sum_{i=1}^q x_i = h_2,  
\end{equation}
 and $J \in \Sym_q$ the sequences $(\mathcal{E}^{\mathbf{x}}(W_{n},J))$ converge, for fixed $q$ and $J$ uniformly in $\mathbf{x}$.

Fix $q$.
Our aim is to find for all $\na \in A_{c_1(q)}$ an $\mathbf{x}$, so that the condition (\ref{feltt}) is satisfied, $\na \in A_\mathbf{x}$ and $A_\mathbf{x} \subseteq A_{c_1(q)}$,
where $c_1(q) \leq x_{i} \leq a_{i}$ for  $i=1, \ldots, q$. As $h_1 < h_2 \leq 1$, there exists such an $\mathbf{x}$ for all $\na \in  A_{c_1(q)}$, let us denote it by $\mathbf{x}_\na$, for convenience set $(\mathbf{x}_\na)_i=\frac{h_1}{q} +\frac{a_i-h_1}{1-h_1}(h_2-h_1)$.
According to this correspondence we have $A_{c_1(q)}= \bigcup_{\na \in  A_{c_1(q)}} A_{\mathbf{x}_\na}$.
So for an arbitrary graphon $W$ and $J \in \Sym_q$ we have
\[
\mathcal{E}^{c_1(q)}(W,J) = \inf_{\na \in  A_{c_1(q)}} \mathcal{E}^{\mathbf{x}_\na}(W,J).
\]
We fix $\varepsilon > 0$, $J \in \Sym_q$, and apply Lemma \ref{ekviv2} for the case that the conditions of the theorem are satisfied. Then there exists a  $n_0 \in \mathbb{N}$, so that for all $n,m > n_0$, for all $\mathbf{x}$ which satisfies (\ref{feltt}), and implies
\[
\left|\mathcal{E}^{\mathbf{x}}(W_n,J) -\mathcal{E}^{\mathbf{x}}(W_m,J)\right|< \varepsilon.
\]
Let $\varepsilon ' > 0$ be arbitrary and $\nb \in A_{c_1(q)}$ such that
$\mathcal{E}^{c_1(q)}(W_m,J) + \varepsilon ' > \mathcal{E}^{\mathbf{x}_\nb}(W_m,J)$. 
Then
\begin{align*}
\mathcal{E}^{c_1(q)}(W_n,J) -\mathcal{E}^{c_1(q)}(W_m,J) &< \mathcal{E}^{c_1(q)}(W_n,J) - \mathcal{E}^{\mathbf{x}_\nb}(W_m,J) + \varepsilon ' \nonumber \\ &\leq
 \mathcal{E}^{\mathbf{x}_\nb}(W_n,J) - \mathcal{E}^{\mathbf{x}_\nb}(W_m,J) + \varepsilon ' <
\varepsilon + \varepsilon '. \nonumber
\end{align*}
The lower bound of $\mathcal{E}^{c_1(q)}(W_n,J) -\mathcal{E}^{c_1(q)}(W_m,J)$ can be established completely similarly and as $\varepsilon '$ was arbitrary, it follows that
\[
\left|\mathcal{E}^{c_1(q)}(W_n,J) -\mathcal{E}^{c_1(q)}(W_m,J)\right|< \varepsilon,
\]
which verifies the statement of the theorem.
\qed

A direct consequence is the version of Theorem \ref{implik} for weighted graphs.

\begin{corollary}
\label{implik2}
Let $G_n$ be a sequence of weighted graphs with uniformly bounded edge weights, and  $\frac{\alpha_{\max}(G_n)}{\alpha_{G_n}} \to 0$ ($n \to \infty$). Let $c_1$ and $c_2$ be two lower threshold functions, so that $c_1(q)q = h_1 < h_2 =c_2(q)q$ for all $q \geq 1$ for some $0 \leq h_1,h_2 \leq 1$.
If for every $q \geq 1$ and $J \in \Sym_q$ the sequences $(\hat{\mathcal{E}}^{c_2(q)}(G_{n},J))_{n \geq 1}$ converge, then for every $q \geq 1$ and $J \in \Sym_q$ the sequences $(\hat{\mathcal{E}}^{c_1(q)}(G_{n},J))_{n \geq 1}$ also converge.
\end{corollary}
The proof of Corollary \ref{implik2} can be easily given through the combination of the results of Theorem \ref{mikrokozel} and Theorem \ref{implik}.

Concluding this section we would like to mention a natural variant of the LTGSEs, the upper threshold ground state energies (UTGSE). Here we will only give an informal description of the definition and the results and leave the details to the reader, everything carries through analogously to the above. 
The homogeneous UTGSE, denoted by $\hat{\mathcal{E}}^{c\uparrow}(G,J)$, is determined by a formula similar to (\ref{defltgse1}) with the set $A_c$ replaced by $A^c$, that is the set of probability distributions whose components are at most $c$, the general variant of the UTGSE is defined in the same manner. The equivalence corresponding to the one stated in Lemma \ref{ekviv2} between the general and the homogeneous version's convergence follows by the same blow-up trick as there, here for $c(q)q = h \geq 1$. The counterpart of Theorem \ref{implik} also holds true in the following form for $1\leq c_2(q)q \leq c_1(q)q\leq q$: If for every $q \geq 1$ and $J \in \Sym_q$ the sequences $(\mathcal{E}^{c_2(q)\uparrow}(W_{n},J))_{n \geq 1}$ converge, then for every $q \geq 1$ and $J \in \Sym_q$ the sequences $(\mathcal{E}^{c_1(q)\uparrow}(W_{n},J))_{n \geq 1}$ also converge. This conclusion comes not unexpected, it says, as in the LTGSE case, that less restriction on the set $A^c$ weakens the convergence property of a graph sequence.  
\section{Counterexamples}
In this section we provide an example of 
 a graphon family whose elements can be distinguished for a larger $c_2(q)$ lower
 threshold function for some pair of $q_0 \geq 1$ and $J_0 \in \Sym_q$ by
 looking at $\mathcal{E}^{c_2(q_0)}(W,J_0)$, but whose LTGSEs are identical for
 some smaller $c_1(q)$ lower threshold function for all $q \geq 1$ and $J \in \Sym_q$.
 Based on this it is possible to construct
 a sequence of graphs, whose $c_1(q)$-LTGES's converge for
 every $q \geq 1$ and $J \in \Sym_q$, but not 
the $c_2(q)$-LTGSEs through the same randomized method presented in \cite{lovlac2} to show a non-convergent graph sequence with 
convergent ground state energies.

In the second part of the section we demonstrate that there exist a family of graphons, where elements can be distinguished from each other by looking only at their LTGSEs for an arbitrary small, but positive $c(q)$ lower threshold function, but whose corresponding GSEs without any threshold are identical. 
\begin{example}\rm{
An example which can be treated relatively easily are block-diagonal graphons which are defined for the parameters $0 \leq \alpha \leq 1$, $0 \leq \beta_1 , \beta_2$ as
\[ W(x,y) = \left\{ \begin{array}{ll}
\beta_1 & \textrm{, if $0 \leq x,y \leq \alpha$}\\
\beta_2 & \textrm{, if $\alpha < x,y \leq 1$}\\
0 & \textrm{, else.}
\end{array} \right. \]

In the case of  $c(q)q=h$ , $1-\alpha \geq h$, $\beta_2=0$, for arbitrary $q \geq 1$ and $J \in \Sym_q$ we have $\mathcal{E}(W,J)=\mathcal{E}^{c(q)}(W,J)$. Choosing $\beta_1=\frac{1}{\alpha^2}$, we get a one parameter family of graphons which have identical $c(q)$-LTGSEs parametrized by $\alpha$ with  $0 < \alpha \leq 1-h$. This means that $\mathcal{E}^{c(q)}(W(\alpha),J) = \mathcal{E}(\mathbb{I},J)$, where $\mathbb{I}$ stands for the constant $1$ graphon.

For every $\alpha_0 > 1-h$ there are $q \geq 1$ and $J \in \Sym_q$, so that the former equality does not hold anymore. Let $J_q \in \Sym_q$ be the  $q \times q$ matrix, whose diagonal entries are $0$,  all other entries being $-1$ (this is the $q$-partition mincut problem). Then  $\mathcal{E}(\mathbb{I},J_q)=0$ for all $q \geq 1$, but for $q_0$ large enough 
$\mathcal{E}^{c(q_0)}(W(\alpha_0),J_{q_0}) >  0$, we leave the details to the reader.
}
\end{example}
With the aid of the previous example it is possible to construct a 
sequence
 of graphs which verify  that in Theorem \ref{implik} the implication of the convergence property of the sequence is strictly one-way.
 This example is degenerate in the sense that the graphs consist of
 a quasi-random part and a sub-dense part with the bipartite graph spanned between the two parts also being sub-dense.

\begin{example}\rm{
Let us consider block-diagonal graphons with $0< \alpha < 1$, $\beta_1,\beta_2 > 0$. It was shown in \cite{lovlac2} that if we restrict our attention to a subfamily of block-diagonal graphons, where  $\alpha^2\beta_1 + (1- \alpha)^2\beta_2$ is constant, then in these subfamilies the corresponding GSEs are identical.  Let $c(q)$ be an arbitrarily small positive threshold function. Next we will show that the $c(q)$-LTGSEs determine the parameters of the block-diagonal graphon at least for a one-parameter family (up to graphon equivalence, since $(\alpha,\beta_1, \beta_2)$ belongs to the same equivalence class as $(1-\alpha,\beta_2,\beta_1)$). The constant $\delta_{ij}$ is $1$, when $i=j$, and $0$ otherwise.

The value of the expression $\alpha^2\beta_1 + (1- \alpha)^2\beta_2$ is determined by the MAXCUT problem by $\Ec(W,J)$ with $q=2$ and $J_{ij} =1- \delta_{ij}$.

In the second step let  $q_0$ be as large so that $c(q_0) < \min(\alpha, 1-\alpha)$ holds, and let $J$ be the $q_0 \times q_0$ matrix with entries $J_{ij}=- \delta_{i1} \delta_{j1}$. In this case $\mathcal{E}(W,J)=0$, but simple calculus gives $-\mathcal{E}^{c(q_0)}(W,J)= - \frac{\beta_1 \beta_2}{\beta_1 + \beta_2} c(q_0)^2$. Hence $\frac{1}{\beta_1} + \frac{1}{\beta_2}$ is determined by the LTGSEs.

The extraction of a third dependency of the parameters from $c(q)$-LTGSEs requires little more effort, we will only sketch details here. First consider $\alpha$'s with $\min(\alpha, 1-\alpha) \geq  c(2)$. Let for $q=2$ and $k \geq 1$ be
\[ 
J_k = \left( \begin{array}{cc}
1 & -k \\
-k & 2 \\
\end{array} \right).
\]
For every $\alpha$ with $\min(\alpha, 1-\alpha) \geq  c(2)$ we have
\[
\lim_{k \to \infty} -\mathcal{E}^{c(2)}(W,J_k) = \alpha^2\beta_1 + (1-\alpha)^2\beta_2 + \max(\alpha^2\beta_1, (1-\alpha)^2\beta_2).
\]

Now apply the notion of the general lower threshold: for $q=2$ let $c_1(n)=2c(2)/n$ and $c_2(n)=2c(2)(n-1)/n$ two threshold functions , let us first  consider the threshold $\mathbf{x}_n =(c_1(n),c_2(n))$. If $\alpha \geq c_1(n)$ or $1- \alpha \geq c_1(n)$, then analogously to the case of the homogeneous lower thresholds
\[
\lim_{k \to \infty} -\mathcal{E}^{\mathbf{x}_n}(W,J_k) = 2\alpha^2\beta_1 + (1-\alpha)^2\beta_2 \quad \textrm{or} \quad \alpha^2\beta_1 + 2(1-\alpha)^2\beta_2.
\]

If for example  $\alpha < c_1(n)$, then it is easy to see that the LTGSE is going to infinity, because for some $k_0$, for all $k>k_0$:
\[
-\mathcal{E}^{\mathbf{x}_n}(W,J_k) < -k (c_1(n) - \alpha) c_2(n) \min\{\beta_1,\beta_2\} + 2(\alpha^2\beta_1 + (1-\alpha)^2\beta_2).
\]
So for  fixed $n$ then 
\[
\lim_{k \to \infty}  -\mathcal{E}^{\mathbf{x}_n}(W,J_k) = - \infty.
\]
To actually be able to extract the expression $\alpha^2\beta_1 + (1-\alpha)^2\beta_2 + \max(\alpha^2\beta_1, (1-\alpha)^2\beta_2)$, we only have to consider the lower threshold obtained by swapping the bounds, $\mathbf{x}_n'=(c_2(n), c_1(n))$.

Then, if $\alpha \geq c_1(n)$ or $1- \alpha \geq c_1(n)$, we have
\begin{align*}
&\max\{\lim_{k \to \infty} -\mathcal{E}^{\mathbf{x}_n}(W,J_k),\lim_{k \to \infty} -\mathcal{E}^{\mathbf{x'}_n}(W,J_k)\} \\ & \qquad =\alpha^2\beta_1 + (1-\alpha)^2\beta_2 + \max\{\alpha^2\beta_1, (1-\alpha)^2\beta_2\},
\end{align*}
otherwise
\[
\max\{\lim_{k \to \infty} -\mathcal{E}^{\mathbf{x}_n}(W,J_k),\lim_{k \to \infty} -\mathcal{E}^{\mathbf{x'}_n}(W,J_k)\}=-\infty
\]
For every $\alpha$ there is a minimal $n_0$ so that one of the conditions $\alpha \geq c_1(n)$ and $1- \alpha \geq c_1(n)$ is satisfied, and for $n < n_0$ the LTGSEs corresponding to $\mathbf{x}_n$ and $\mathbf{x'}_n$ tend to infinity when $k$ goes to infinity. Therefore the expression $\alpha^2\beta_1 + (1-\alpha)^2\beta_2 + \max(\alpha^2\beta_1, (1-\alpha)^2\beta_2)$ is determined by $c(q)$-LTGSEs.

Consider the one-parameter block-diagonal graphon family analyzed in \cite{lovlac2}, that is $W(\alpha)=W(\alpha, \frac{1}{\alpha},\frac{1}{1-\alpha})$, where $0< \alpha <1$. In this case the values of our first two expressions are constant, for every $0< \alpha <1$ we have
\begin{align*}
\frac{1}{\beta_1} + \frac{1}{\beta_2} &= 1, \\
\alpha^2 \beta_1 + (1-\alpha)^2 \beta_2 &= 1.
\end{align*}
But by applying the third expression for  $c(q)$-LTGSEs, we extract $\max\{\alpha^2\beta_1,(1-\alpha)^2 \beta_2\}=\max\{\alpha,1-\alpha\}$, which determines the graphon uniquely in this family up to equivalence.}
\end{example}

\end{document}